# Onsite Job Scheduling by Adaptive Genetic Algorithm

Avijit Basak          Subhas Acharya

## ABSTRACT

Onsite Job Scheduling is a specialized variant of Vehicle Routing Problem (VRP) with multiple depots. The objective of this problem is to execute jobs requested by customers, belonging to different geographic locations by a limited number of technicians, with minimum travel and overtime of technicians. Each job is expected to be completed within a specified time limit according to the service level agreement with customers. Each technician is assumed to start from a base location, serve several customers and return to the starting place. Technicians are allotted jobs based on their skill sets, expertise levels of each skill and availability slots. Although there are considerable number of literatures on VRP we do not see any explicit work related to Onsite Job Scheduling. In this paper we have proposed an *Adaptive Genetic Algorithm* to solve the scheduling problem. We found an optimized travel route for a substantial number of jobs and technicians, minimizing travel distance, overtime duration as well as meeting constraints related to SLA.

## 1. Introduction

Vehicle Routing Problem is among the most widely studied problem in the field of operational research. VRP also attracts special interests from researchers due to its environmental and economic importance. VRP attempts to solve delivery problem where goods are delivered from few specified depots to designated customers by vehicles. It attempts to minimize the travelled route, number of vehicles needed, meeting agreed SLA with customers, and maximize collection/profit. There are different variants of VRP like VRP with pickup and delivery, VRP with time windows, capacitated VRP, VRP with multiple trips, multi-depot VRP, VRP with transfer. We have studied an onsite job scheduling problem which is a special variant of VRP, where each technician visits job premises in an assigned order, executes the jobs and returns to their starting location. We have tried to minimize the total travelled distance, overtime duration while meeting all job SLAs. Historically diverse types of algorithms have been used to solve VRP problems as described in the references below. We have chosen genetic algorithm to find a near optimal solution to this NP-hard constrained optimization problem.

Genetic algorithm [1] was originally proposed by John Holland following the Survival of Fittest principle of Sir Charles Darwin. The algorithm works based on three operators i.e., selection, crossover, and mutation on a population of randomly generated chromosomes and eventual convergence of population over generations. However, crossover and mutation based on fixed probability does not provide good result. In this paper, we have used rank based adaptive approach[5] for crossover and mutation probability generation to improve quality of optimization.

## 2. Problem Description

In this onsite job scheduling problem, we have addressed the problem of execution of jobs for individual customers at scattered locations by a fixed number of technicians in the city of Ahmedabad, India. We have conceptualized few requirements and constraints as part of the problem. The worker to job ratio is maintained as 1:4. We have not assumed any existence of depot. Each job is considered to have certain priority and an expected completion duration. We have considered ten different skills for all jobs. Each job and technician can have one or two skill sets. Technicians are also assumed to have their availability time slots. Technicians were considered to travel from their registered base location to the assigned job locations according to the scheduled job sequence, execute the jobs and return to their base location. Considering both single and multi-skill jobs, each job is allotted to one and only one worker based on required skill set and expertise level of the workers. A buffer has been considered as a multiplicative factor of expected job completion duration based on worker skill level. The job execution must be completed following the customer SLA constraints to avoid any penalty. Each technician is assumed to travel independently using their own travel arrangements. Our primary objective was to optimize the total amount of distance the technicians travel each day along with overtime duration and at the same time meet the SLA constraints. We have paid maximum priority on executing jobs within predetermined SLA.

## 3. Solution Approach

Onsite scheduling problem is a muti-objective problem with an objective to minimize travelled distance, job tardiness cost and overtime of workers. There are different approaches to solve a multi-objective problems like **a priori method** and **a posteriori method**. A priori method requires sufficient preference information before the start of evaluation process. A posteriori method aims to obtain solutions representing the Pareto front. In this solution we have studied the scenario of automated scheduling with no decision maker hence assumed enough information should have been available to determine the preferences of each evaluation function before start of optimization. So, we adopted linear scalarization approach as **a priori method** to evaluate the solution.

We have adopted random key encoding to represent job sequence as chromosome. Workers are assigned individual jobs based on job skill set requirement. Worker allocation is maintained as a mapping between decoded job id and worker id. The primary approach we have adopted to optimize the schedule is optimization of travel distance for each worker by crossover and alteration of worker allocation to individual jobs by mutation. The entire optimization process has been described by the flowchart in Figure 1. For a pictorial description of the detailed process, we have considered six jobs with two skill sets and three workers as shown in the tables (Figure 2).

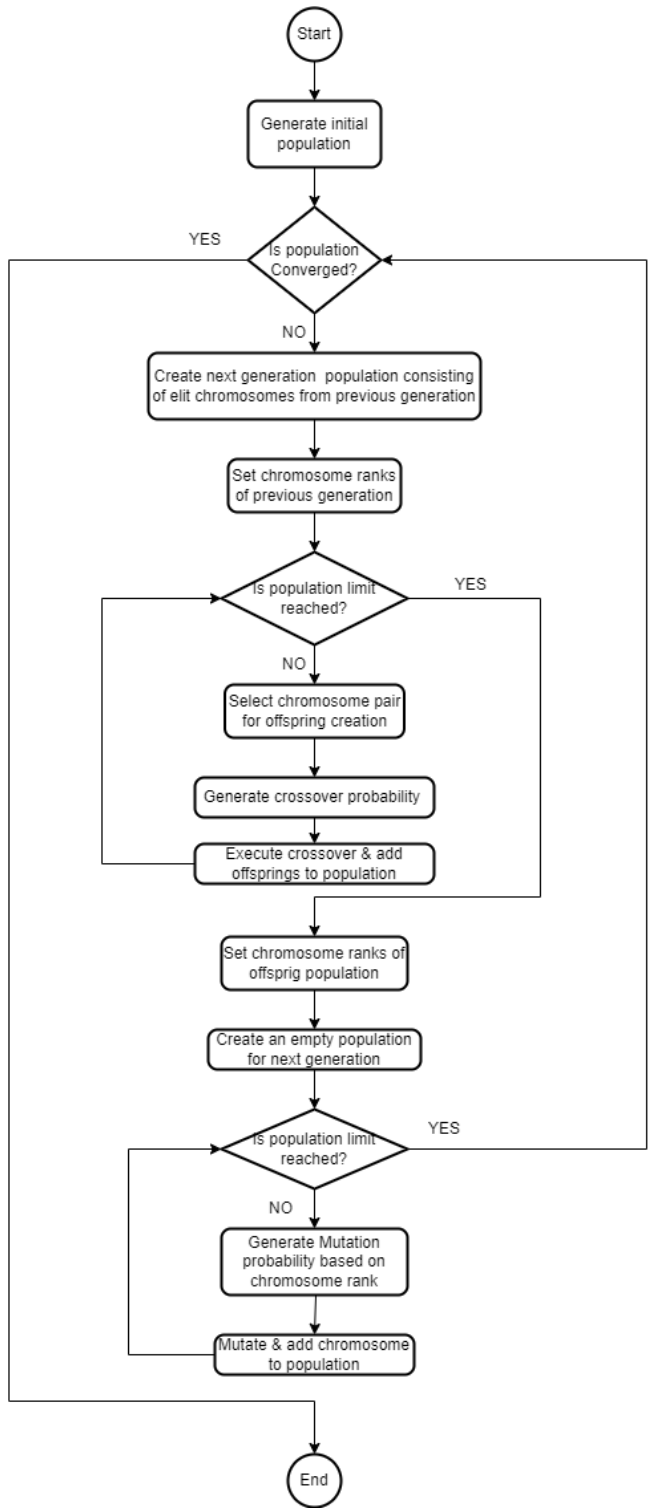

*Figure 1*

| Job to Skillset mapping | |
|---|---|
| Job Id | Skill Set Id |
| 1 | 1 |
| 2 | 1 |
| 3 | 1 |
| 4 | 1 |
| 5 | 2 |
| 6 | 2 |

| Worker to Skillset mapping | |
|---|---|
| Worker Id | Skill Set Id |
| 1 | 1 |
| 2 | 1 |
| 3 | 2 |

*Figure 2*

### 3.1. Schedule Encoding

We have encoded schedules as individual chromosomes. Chromosomes are designed to contain information of job sequence and allocated workers to the jobs. The job sequence is represented using a random key encoding. A set of random numbers (between 0-1) are chosen to represent job Ids as part of initial population generation. The workers are chosen randomly based on the job skill requirement and are kept using a map of job id and worker id. A sample chromosome structure is presented below in encoded and decoded form.

**Encoded View**

| Schedule Chromosome | .3 | .7 | .2 | .33 | .99 | .65 |
|---|---|---|---|---|---|---|
| | 1 | 3 | 2 | 1 | 3 | 2 |

**Decoded View**

| Job Id | 2 | 5 | 1 | 3 | 6 | 4 |
|---|---|---|---|---|---|---|
| Worker Id | 1 | 3 | 2 | 1 | 3 | 2 |

*Figure 3*

### 3.2. Selection

Different selection mechanisms have been proposed for genetic algorithm like fitness proportionate selection, rank based selection, tournament selection[7] etc. Among them tournament selection shows superior performance and optimization quality in multimodal problem domain. In this work we have used tournament selection with a tournament size of 10% population size.

### 3.3. Elitism

Elitism was originally proposed by De Jong in his seminal dissertation[6]. Elitism ensures the best chromosomes are not lost during crossover over generations hence improves the optimization quality considerably. We have used elitism with 10% elitism percentage.

## 3.4. Crossover

We have used one point crossover as part of the optimization process due to its simplistic nature. The crossover process is used to optimize the job execution sequence of each worker. Worker to job allocation mapping is kept unaltered in the offspring chromosomes. The diagram below describes the crossover process in both encoded and decoded views. The crossover probability is determined using a rank based adaptive approach based on maximum rank of mating chromosomes as shown in the equation below.

$$p_c = p_{c-MIN} + (p_{c-MAX} - p_{c-MIN}) * (1 - \frac{\max(r_{ch-1}, r_{ch-2}) - 1}{N-1})$$

Where  $p_{c-MAX}$ = Maximum crossover probability

$p_{c-MIN}$ = Minimum crossover probability

$r_{ch-1}$ = Rank of chromosome 1

$r_{ch-2}$ = Rank of chromosome 2

$\max(r_{ch-1}, r_{ch-2})$ = Maximum of $r_{ch-1}$ and $r_{ch-2}$

$N$ = No of chromosome in population

$p_c$ = Crossover probability

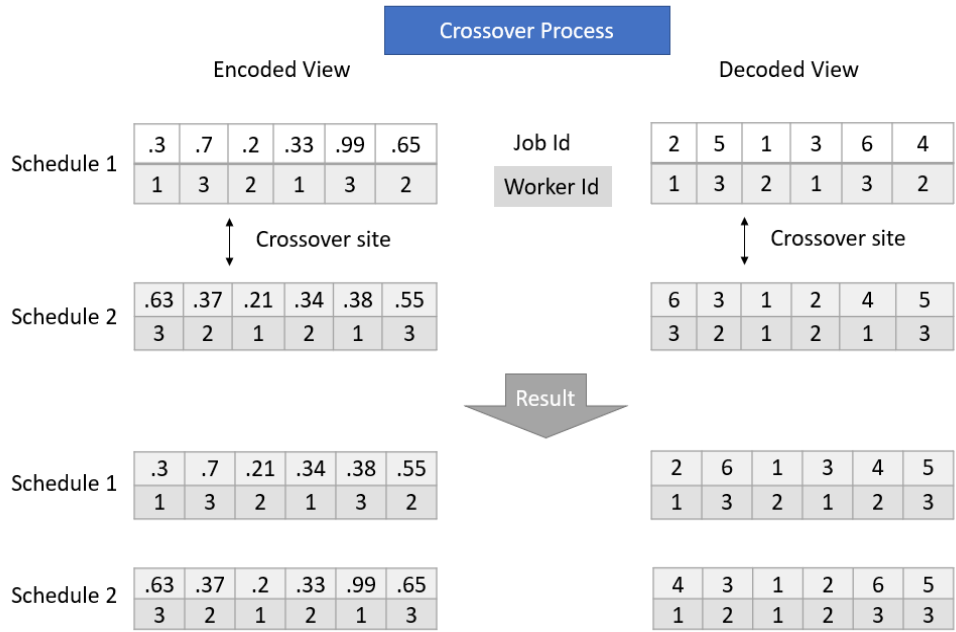

Figure 4

## 3.5. Mutation

We have used mutation to explore better allocation of worker to jobs. As part of the process, we have selected few random jobs based on mutation probability. Workers are allocated randomly depending on the required job skills for those selected jobs. This process helped us to explore different sort of worker allocations and hence alter travel cost as well as job tardiness cost. Mutation probability was determined using a rank based adaptive approach following the equation below.

$$p_m = p_{m-MIN} + (p_{m-MAX} - p_{m-MIN}) * (1 - \frac{r_{ch}-1}{N-1})$$

Where  $p_{m-MAX}$ = Maximum mutation probability

$p_{m-MIN}$ = Minimum mutation probability

$r_{ch}$ = Rank of chromosome

$N$ = No of chromosome in population

$p_m$ = Mutation probability

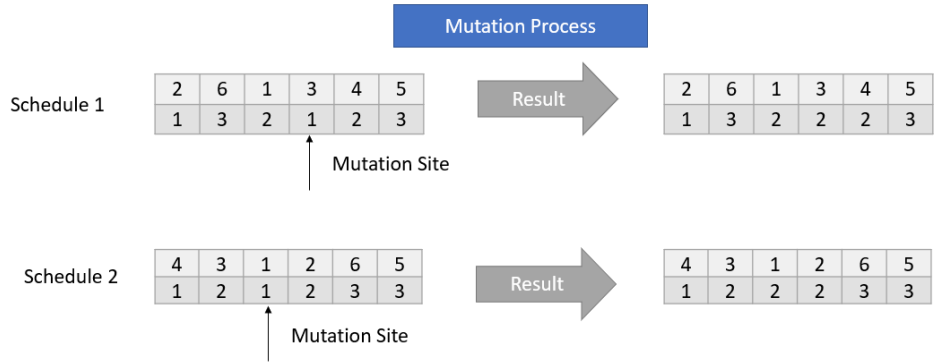

Figure 5

## 3.6. Cost Function

Weighted sum of normalized travel distance, job completion delay time and overtime are used to calculate the cost function of a schedule in this application. Job completion delay time is calculated as difference between expected job completion time and job SLA. Cost due to job delay time is calculated as exponential of job completion delay in minutes. Schedules having positive delay for any job is rejected as part of the process. Overtime of each worker is calculated as a difference between worker's actual working time and contracted regular working time. A maximum value is decided for each component to calculate the normalized values.

$$f_c = w_d * \sum_{i=1}^{m} \frac{d_i}{d_{MAX}} + w_{SLA} \sum_{i=1}^{n} \frac{P_i}{P_{avg}} * e^{(t_i-T_i)/T_{MAX}} + w_T * \sum_{i=1}^{m} \frac{O_i}{O_{MAX}}$$

where,  $m$ = No of workers

$n$ = No of Jobs

$d_i$ = Distance travelled by i$_{th}$ worker

$d_{MAX}$ = Maximum distance allowed to be travelled by a single worker

$t_i$ = Expected completion time of i$_{th}$ Job

$T_i$ = SLA of i$_{th}$ Job

$T_{MAX}$ = Max allowed SLA for any Job

$P_i$ = Priority of i$_{th}$ Job

$P_{avg}$ = Average priority of the Jobs

$w_d$ = Weight of normalized distance travelled

$w_{SLA}$ = Weight of normalized delay cost

$w_T$ = Weight of normalized overtime cost

$O_i$ = Overtime of ith worker

$O_{MAX}$ = Maximum overtime of any worker

## 4. Simulation Result

For simulation purpose we have generated all worker and job locations (by latitude and longitude) randomly within city of Ahmedabad, India. We have considered ten different skills for all jobs and allocated one/two skills randomly to each technician. Skill level is also chosen randomly between five to ten. We have chosen a fixed availability slot from 9 A.M. to 8 P.M for each technician. Each job is assigned priority ranging from one to ten. We have also assumed an expected job completion duration ranging between ten to sixty minutes. We have assumed a multiplicative factor of 20% considering the worker skill level. Other algorithm parameters are listed below. We have provided graphical simulation results for four sets of technicians and jobs. It is evident from the graph plot that for each scenario the algorithm was able to achieve proper optimization.

### 4.1. Algorithm Parameter Values

Maximum distance allowed to be travelled by a single worker ($d_{MAX}$) = 100 km

Max allowed SLA for any Job ($T_{MAX}$) = 1440 mins

Max allowed overtime for any worker ($T_{MAX}$) = 1440 mins

Max allowed overtime for any worker ($O_{MAX}$) = 120 mins

Average priority of Jobs ($P_{avg}$) = 5

Weight of normalized travel cost ($w_d$) = .5

Weight of normalized SLA cost ($w_{SLA}$) = .3

Weight of normalized overtime cost ($w_T$) = .2

Minimum mutation probability ($p_{m-MIN}$) = 0.0

Maximum mutation probability ($p_{m-MAX}$) = 0.2

Elitism Rate= 0.1

Tournament Size= 10% of population size

### 4.2. Scenario 1

Number of Jobs (n) = 80

Number of Workers (m) = 20

Population size (N)= 100

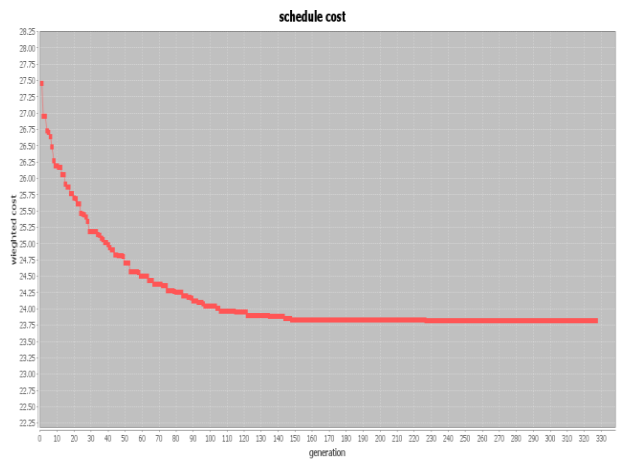
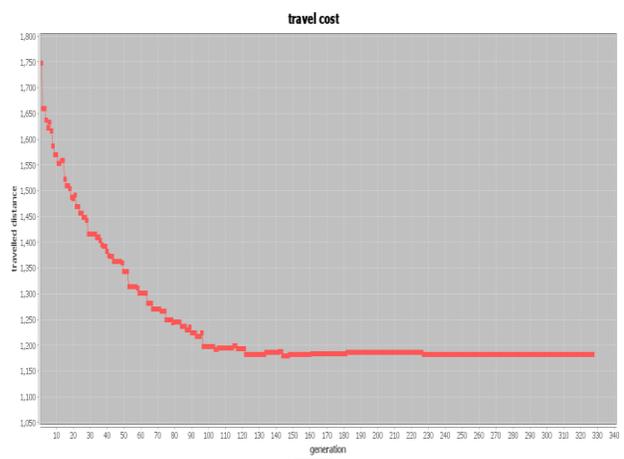
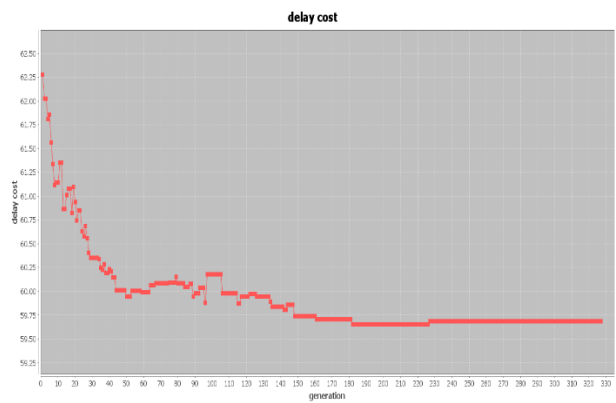
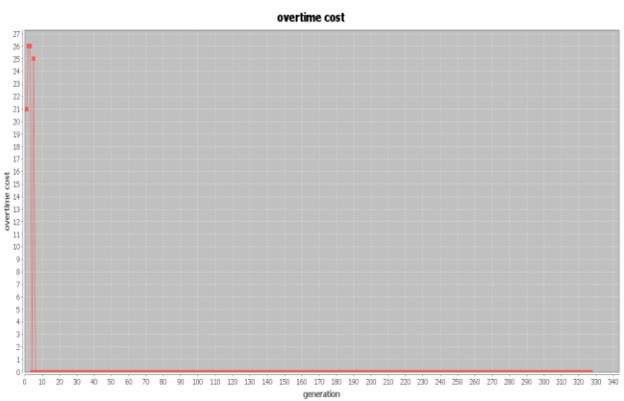

### 4.3. Scenario 2:

Number of Jobs (n) = 160

Number of Workers (m) = 40

Population size (N)= 200

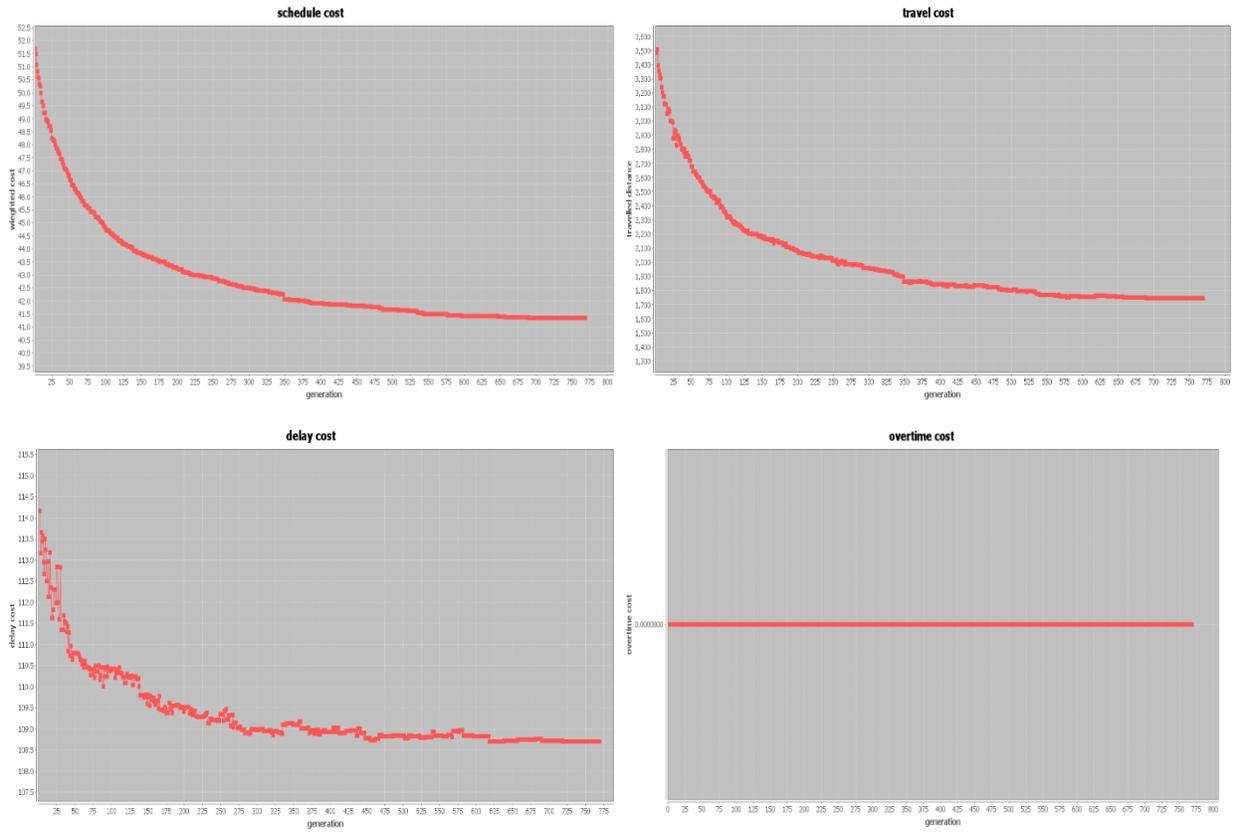

### 4.4. Scenario 3:

Number of Jobs (n) = 320

Number of Workers (m) = 80

Population size (N)= 400

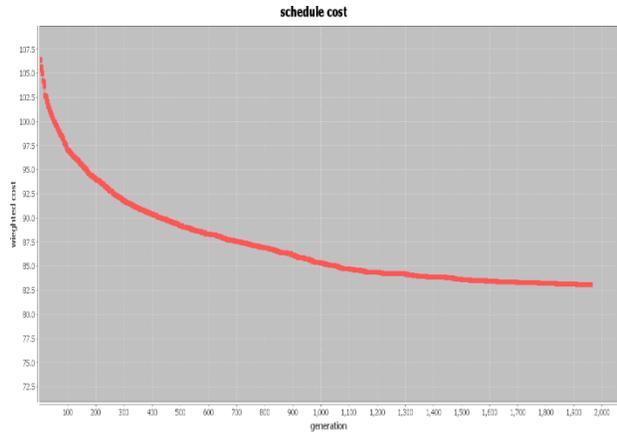

schedule cost

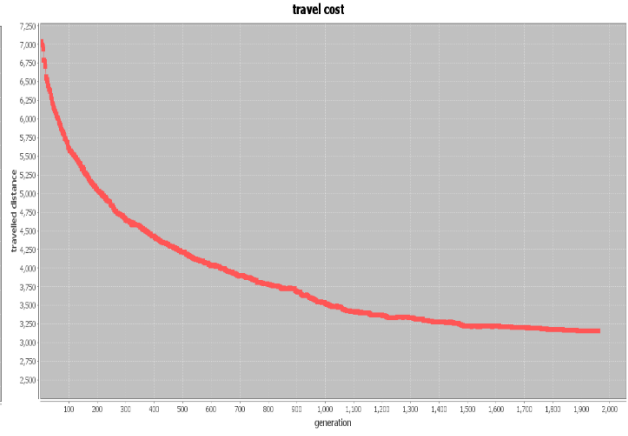

travel cost

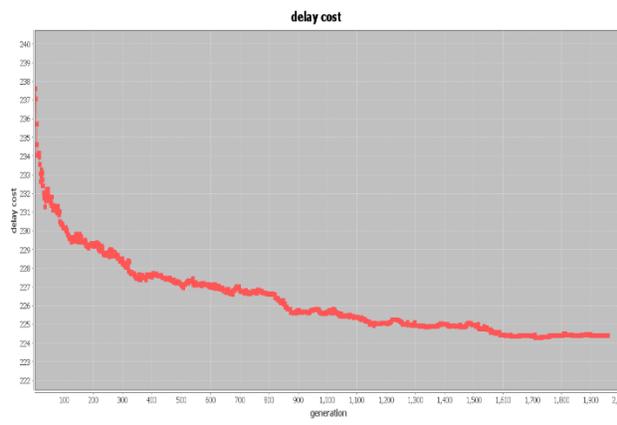

delay cost

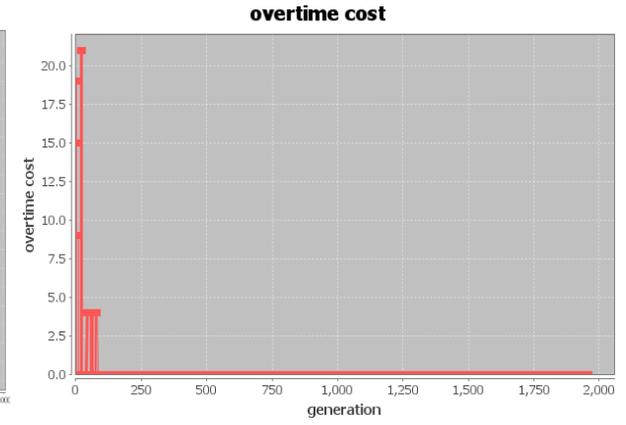

overtime cost

### 4.5. Scenario 4:

Number of Jobs (n) = 400

Number of Workers (m) = 100

Population size (N)= 500

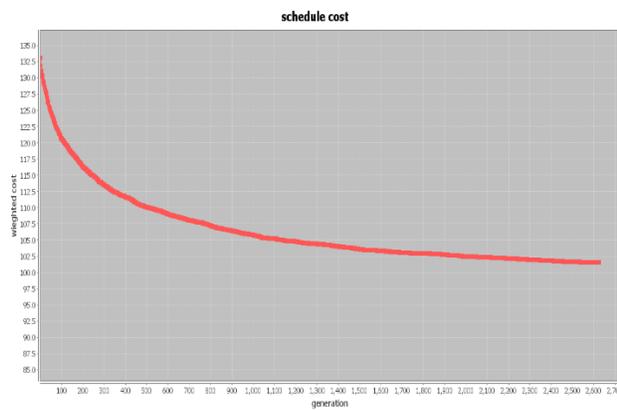

schedule cost

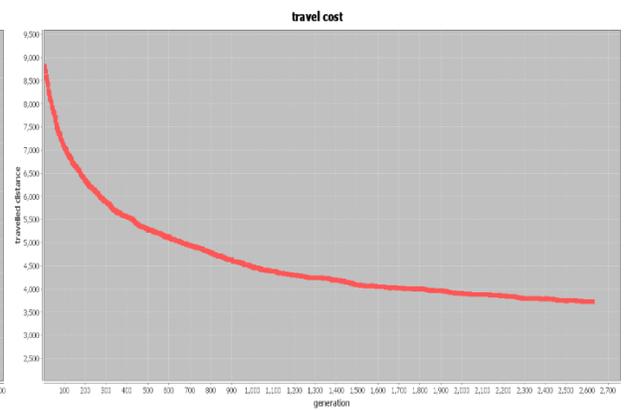

travel cost

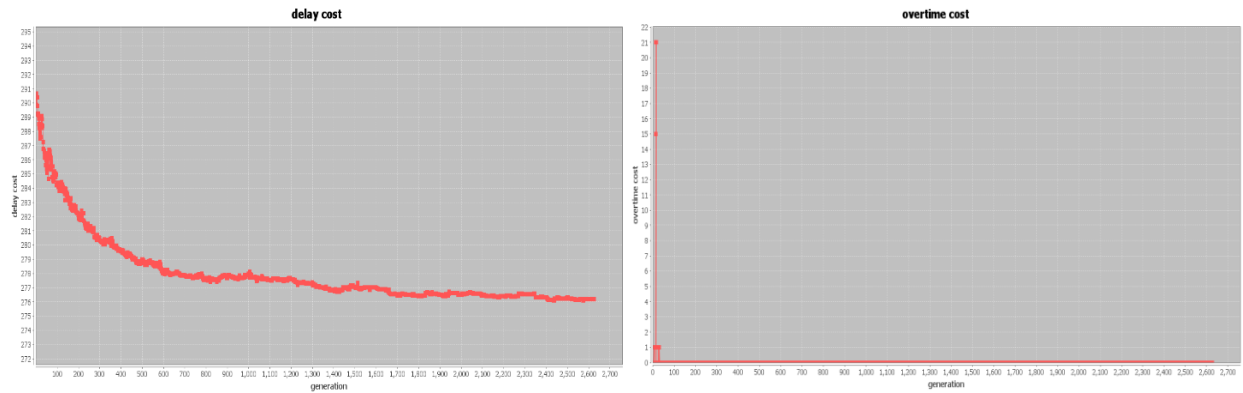

## 5. Conclusion

In this paper we have attempted to solve the onsite scheduling problem for a bulk number of jobs. Although there are many applications of metaheuristics in solving different types of transportation and scheduling problem based on algorithms such as – Genetic Algorithm[1], Particle swarm optimization[4], Adaptive large neighborhood search[10] etc., this paper tries to focus only on adaptive Genetic Algorithm. Rank based adaptive approach for both crossover and mutation probability generation has been adopted to improve the quality of the optimization result. The solution could be utilized by any organization, to create an onsite schedule plan bounded by a specific geographic area automatically, for their resource planning and automated decision making process.

This work can be further extended to compare the results between different heuristics-based approaches over a common set of problem. This approach could also be applied to other specialized areas of VRP like VRP with pickup and delivery, capacitated VRP, VRP with multiple trips, VRP with transfer etc. along with real-time scheduling.